# NONLINEAR STOCHASTIC WAVE EQUATIONS: BLOW-UP OF SECOND MOMENTS IN $L^2$-NORM

BY PAO-LIU CHOW[1]

*Wayne State University*

The paper is concerned with the problem of explosive solutions for a class of nonlinear stochastic wave equations in a domain $\mathcal{D} \subset \mathbb{R}^d$ for $d \leq 3$. Under appropriate conditions on the initial data, the nonlinear term and the noise intensity is proved in Theorem 3.1 that the $L^2$-norm of the solution will blow up at a finite time in the mean-square sense. An example is given to show an application of the theorem.

**1. Introduction.** Consider the Cauchy problem for a nonlinear wave equation

$$(1.1) \qquad \begin{cases} \partial_t^2 u(x,t) = \nabla^2 u + f(u), & t > 0, \\ u(x,0) = g(x), \quad \partial_t u(x,0) = h(x), & x \in \mathbb{R}^d, \end{cases}$$

where $\partial_t = \frac{\partial}{\partial t}$, $\nabla^2$ is the Laplacian operator, and the functions $f, g$ and $h$ are given such that the Cauchy problem (1.1) has a unique local solution. It was first shown by Keller [6] in 1957 that, for a certain class of nonlinear functions $f(u)$, the solution of equation (1.1) becomes infinite or explodes at a finite time, provided that the initial data satisfies appropriate conditions. His result was later generalized by Glassey [4] and others. Since then, it has become known that solutions to more general nonlinear hyperbolic equations may develop singularities in finite time [5], physically manifested as shock waves, tsunami or explosion. It is therefore of interest to examine the effect of a random perturbation to equation (1.1) on the existence of an explosive solution. To study this type of problem, it is necessary to employ some analytical tools from the theory of stochastic partial differential equations (see, e.g., [3]).

Received December 2008; revised March 2009.

[1]Supported in part by the NSF Grant 0600537.

In our previous papers [1, 2], we considered the local and global solutions of a certain class of nonlinear stochastic wave equations and their asymptotic behavior. As an example of the nonexistence of a global solution, we showed that, under some explicit conditions on the initial data and the noise term, the solution of a cubically nonlinear wave equation perturbed by an additive noise will blow up at a finite time in the mean-square sense. In the present paper, it will be shown that this blowup result can be extended to a general class of nonlinear stochastic wave equations given by (2.2). Related to our work, we should mention an earlier paper of Mueller [7] who investigated the global solution of the stochastic wave equation

$$(1.2) \quad \begin{cases} \partial_t^2 u(x,t) = \nabla^2 u + a(u)\,\partial_t W(x,t), & t>0, \\ u(x,0) = g(x), \quad \partial_t u(x,0) = h(x), & x \in \mathbb{R}^d \end{cases}$$

for $d = 1, 2$, where $W(x,t)$ is a Wiener random field. He found a growth condition on $a(u)$, so equation (1.2) has a unique global solution. But the existence of explosive solutions was not considered there. Notice that, in contrast with our equation (2.2), the equation (1.2) is a random perturbation of a linear wave equation. However, its noise amplitude $a(u)$ may be unbounded.

In what follows, we shall first recall some basic results for nonlinear stochastic wave equations in Section 2. Then, in Section 3, we prove the main theorem (Theorem 3.1) which states that, under some sufficient conditions on the initial data and the noise term, the solutions to a class of nonlinear stochastic wave equations whose $L^2$-norms will explode at a finite time in the mean-square sense. Finally, in Section 4, we apply the theorem to a two-dimensional problem to obtain a set of explicit conditions for explosive solutions. In passing, it should be pointed out that we have only obtained the blowup result in a mean-square-$L^2$ sense. The question concerning the existence of explosive solutions, in the almost sure sense and at each point in space, is a challenging problem which remains open.

**2. Preliminaries.** Let $\mathcal{D} \subset \mathbb{R}^d$ be a bounded domain with a smooth boundary $\partial \mathcal{D}$. We set $H = L^2(\mathcal{D})$ with the inner product and norm are denoted by $(\cdot,\cdot)$ and $\|\cdot\|$, respectively. Let $H^1 = H^1(\mathcal{D})$ be the $L^2$-Sobolev space of first order and denote by $H_0^1$ the closure in $H^1$ the space of $C^1$-functions with compact support in $\mathcal{D}$. Let $W(x,t)$ be a continuous Wiener random field defined in a complete probability space $(\Omega, \mathcal{F}, P)$ with a filtration $\mathcal{F}_t$ ([3], p. 38). It has mean $EW(x,t) = 0$ and covariance function defined by

$$EW(x,t)W(y,s) = (t \wedge s) r(x,y), \qquad x, y \in \mathcal{D},$$

where $(t \wedge s) = \min(t,s)$ for $0 \leq t, s \leq T$. In this paper, the spatial correlation function $r(x,y)$ is assumed to be bounded and continuous for $x, y \in D$ and, if $D$ is unbounded,

$$\int_D r(x,x)\,dx < \infty. \tag{2.1}$$

Let $f(\mu)$ and $\sigma(\mu,\xi,x,t)$ be two nonlinear continuous functions for $\mu \in \mathbb{R}$, $\xi \in \mathbb{R}^d, x \in \mathcal{D}$ and $t \geq 0$. We consider the initial-boundary value problem for the nonlinear stochastic wave equation:

$$\begin{cases} \partial_t^2 u(x,t) = (c^2\nabla^2 - \alpha)u + f(u) \\ \qquad\qquad + \sigma(u, Du, x, t)\,\partial_t W(x,t), & t > 0, \\ u(x,0) = g(x), \quad \partial_t u(x,0) = h(x), & x \in \mathcal{D}, \\ u(\cdot,t)|_{\partial\mathcal{D}} = 0, \end{cases} \tag{2.2}$$

with $g \in H^1, h \in H$, where $c$ and $\alpha$ are positive parameters, and $D = \partial_x$ denotes the gradient operator.

To consider (2.2) as an Itô equation in $H$, we set $u_t = u(\cdot,t), v_t = v(\cdot,t)$ and so on, and rewrite it as

$$\begin{cases} du_t = v_t\,dt, \\ dv_t = [(\nabla^2 - \alpha)u_t + f(u_t)]\,dt + dM_t(u), & 0 < t < T, \\ u_0 = g, \quad v_0 = h, \end{cases} \tag{2.3}$$

which can be written as the integral equation

$$\begin{cases} u_t = u_0 + \int_0^t v_s\,ds, \\ v_t = v_0 + \int_0^t (\nabla^2 - \alpha)u_s\,ds + \int_0^t f(u_s)\,ds + M_t(u), \end{cases} \tag{2.4}$$

where we set

$$M_t(u) = \int_0^t \sigma_s(u_s)\,dW_s, \tag{2.5}$$

and $\sigma_t(u) = \sigma(u, Du, \cdot, t)$. The stochastic integral (2.5) exists under the above conditions on $r(x,y)$ (see [1], the Appendix).

Suppose that the nonlinear functions $f(\mu)$ and $\sigma(\mu,\xi,x,t)$ are locally Lipschitz continuous with a polynomial growth. In addition, assume that $d \leq 3$ and there exists an energy bound for equation (2.2):

$$\int_0^t (v_s, f(u_s))\,ds + \int_0^t \mathrm{Tr}\,Q_s(u_s)\,ds \tag{2.6}$$
$$\leq c_1 \mathbf{e}(u_t; v_t) + c_2 \int_0^t \mathbf{e}(u_s; v_s)\,ds + c_3 \quad \text{for } 0 \leq t \leq T,$$

where $c_1, c_2, c_3$ are some constants with $c_1 < 1$; $\mathbf{e}(u;v)$ is the energy function defined by

$$\mathbf{e}(u;v) = c^2\|Du\|^2 + \alpha\|u\|^2 + \|v\|^2 \quad \text{for } u \in H^1, v \in H, \tag{2.7}$$

and

(2.8) $$TrQ_t(u) = \int_{\mathcal{D}} r(x,x)\sigma^2(u, Du, x, t)\, dx.$$

Then, by [2], Theorem 3.1, it has a unique continuous global solution $u_t \in H_0^1$ with $\partial_t u_t \in H$, $t \in [0,T]$, for any $T > 0$. Moreover, the following energy equation holds

(2.9) $$\mathbf{e}(u_t; v_t) = \mathbf{e}(u_0; v_0) + 2\int_0^t (v_s, f(u_s))\, ds + \int_0^t TrQ_s(u_s)\, ds$$
$$+ 2\int_0^t (v_s, \sigma_s(u_s)\, dW_s) \qquad \text{a.s., for } 0 \leq t \leq T.$$

We remark that, in order to allow a polynomial nonlinearity, the restriction $d \leq 3$ is needed to prove the theorem cited above. Now, we define

(2.10) $$\phi(t) = \tfrac{1}{2} E[(u_t, u_t)] = \tfrac{1}{2} E\|u_t\|^2.$$

Then, similar to [1], Lemma 3.1, we can compute its first two derivatives

(2.11) $$\phi'(t) = E[(u_t, v_t)] = (u_0, v_0) + E\int_0^t \{\|v_s\|^2 - c^2\|Du_s\|^2 - \alpha\|u_s\|^2 + (u_s, f(u_s))\}\, ds,$$

(2.12) $$\phi''(t) = E\{\|v_t\|^2 - c^2\|Du_t\|^2 - \alpha\|u_t\|^2 + (u_t, f(u_t))\}.$$

On the other hand, without the energy bound, one can only assert the existence of a unique local solution. That is, there exist a unique solution $u_t$ as before for $t < \tau$ and a sequence of stopping times $\{\tau_n\}$ such that $\tau_n \uparrow \tau < \infty$ and the following holds for any $n$:

(2.13) $$\begin{cases} u_{t \wedge \tau_n} = u_0 + \int_0^{t \wedge \tau_n} v_{s \wedge \tau_n}\, ds, \\ v_{t \wedge \tau_n} = v_0 + \int_0^{t \wedge \tau_n} (\nabla^2 - \alpha) u_{s \wedge \tau_n}\, ds \\ \qquad\qquad + \int_0^{t \wedge \tau_n} f(u_{s \wedge \tau_n})\, ds + M_{t \wedge \tau_n}(u). \end{cases}$$

In particular, the limiting stopping time $\tau$ may be a positive constant. It is clear that if $P\{\tau = \infty\} = 1$, then the local solution becomes a global solution.

**3. Explosive solutions.** Consider the local solution of the problem (2.2) in $\mathcal{D} \in \mathbb{R}^d$ with $d \leq 3$. We shall show that, under a set of sufficient conditions, the second-moment $E\|u_t\|^2$ of the local solution will explode at a finite time in the $L^2$-norm. To this end, we first impose the following conditions:

(A1) The function $f(\mu)$, for $\mu \in \mathbb{R}$, is locally Lipschitz continuous and may be of polynomial growth.

(A2) For $\mu \in \mathbb{R}$, $\xi \in \mathbb{R}^d, x \in \mathcal{D}$ and $t \geq 0$, the function $\sigma(\mu, \xi, x, t)$ is continuous in $x, t$ and it is locally Lipschitz continuous in $\mu$ and $\xi$.
(A3) $\int_0^\infty \int_\mathcal{D} r(x,x) q(x,t) \, dx \, dt < \infty$, where $q(x,t) = \sup_{\mu,\xi} \sigma^2(\mu, \xi, x, t)$, for $\mu \in \mathbb{R}$ and $\xi \in \mathbb{R}^d$.

We will state and prove the following theorem concerning explosive solutions.

THEOREM 3.1. *Let conditions* (A1)–(A3) *hold true. Suppose that, for $u_0 \in H^1, v_0 \in H$, the problem (2.2), or the system (2.3) has a unique continuous local solution $u_t \in H_0^1$ with $\partial_t u_t \in H$. Then there exists a time $T_e > 0$ such that*

$$\lim_{t \to T_e^-} E\|u_t\|^2 = +\infty, \tag{3.1}$$

*provided that the following conditions are satisfied:*

(B1) *The inner product* $(u_0, v_0) > 0$.
(B2) $(F(u_0), 1) > \dfrac{1}{2}\Big\{ \mathbf{e}(u_0; v_0) + \displaystyle\int_0^\infty \int_\mathcal{D} r(x,x) q(x,t) \, dx \, dt \Big\}$.
(B3) $(u, f(u)) \geq \frac{1}{2}(F(u), 1)$, *for any bounded continuous function $u$ on $\mathcal{D}$, where $F(u) = \int_0^u f(\mu) \, d\mu$ and $(F(u), 1) = \int_\mathcal{D} F(u(x)) \, dx$.*

PROOF. We shall prove the theorem by contradiction. Suppose the solution is global so that, for any $T > 0$,

$$E\|u_t\|^2 < \infty \quad \text{for } 0 \leq t \leq T. \tag{3.2}$$

Define $\phi(t) = \frac{1}{2} E\|u_t\|^2$ by (2.10). We will show that $\phi(t)$ becomes unbounded in finite time, or the inequality (3.2) is false. Alternatively, for any $\lambda > 0$, it suffices to show that

$$\psi(t) = \phi^{-\lambda}(t) \tag{3.3}$$

goes to zero in finite time. To proceed we first compute the first two derivatives of $\psi$ to obtain

$$\psi'(t) = -\lambda \phi^{-(\lambda+1)}(t) \phi'(t) = -\lambda \phi^{-(\lambda+1)}(t) E[(u_t, v_t)] \tag{3.4}$$

and

$$\begin{aligned}\psi''(t) &= -\lambda \phi^{-(\lambda+1)}(t) \Big\{ \phi''(t) - (\lambda+1) \frac{[\phi'(t)]^2}{\phi(t)} \Big\} \\ &\leq -\lambda \phi^{-(\lambda+1)}(t) \{ \phi''(t) - 2(\lambda+1) E\|v_t\|^2 \},\end{aligned} \tag{3.5}$$

where use was made of (2.11) and the inequality: $[\phi'(t)]^2 \leq E\|u_t\|^2 E\|v_t\|^2$.

Notice that
$$\int_0^t (v_s, f(u_s))\, ds = \int_{\mathcal{D}} [F(u_t) - F(u_0)]\, dx = (F(u_t) - F(u_0), 1).$$

From this equation together with the energy equations (2.9) and (2.7), we can deduce that

(3.6)
$$E\|v_t\|^2 = \left\{\mathbf{e}(u_0; v_0) - 2(F(u_0), 1) + E\int_0^t TrQ_s(u_s)\, ds\right\}$$
$$+ E\{2(F(u_t), 1) - c^2\|Du\|^2 - \alpha\|u\|^2\} \quad \text{for } 0 \le t \le T.$$

Upon substituting equations (2.12) and (3.6) into the inequality (3.5) and after some simplification, we can obtain

(3.7)
$$\psi''(t) \le -\lambda\phi^{-(\lambda+1)}(t)\left\{2\lambda\left[2(F(u_0),1) - \mathbf{e}(u_0;v_0) - E\int_0^t TrQ_s(u_s)\,ds\right]\right.$$
$$+ (2\lambda - 1)E(c^2\|Du\|^2 + \alpha\|u\|^2)$$
$$\left. + E[(u_t, f(u_t)) - 4\lambda(F(u_t), 1)]\right\}.$$

By choosing $\lambda = 1/2$ in (3.7) and making use of condition (A3), we get

(3.8)
$$\psi''(t) \le -\frac{1}{2}\phi^{-3/2}(t)\left\{\left[2(F(u_0), 1) - \mathbf{e}(u_0; v_0)\right.\right.$$
$$\left. - \int_0^\infty \int_{\mathcal{D}} r(x,x) q(x,t)\, dx\, dt\right]$$
$$\left. + E[(u_t, f(u_t)) - 2(F(u_t), 1)]\right\},$$

which, in view of conditions (B2) and (B3), implies that

(3.9) $$\psi''(t) > 0 \quad \text{for } 0 \le t \le T.$$

In the meantime, by setting $\lambda = 1/2$ in (3.3) and (3.4) and taking (2.10) and (2.11) into account, we can obtain

(3.10) $$\psi(0) = \frac{\sqrt{2}}{\|u_0\|} > 0$$

and

(3.11) $$\psi'(0) = -\frac{\sqrt{2}}{\|u_0\|^3}(u_0, v_0) < 0,$$

by invoking condition (B1). Take $T > T_0 = -\psi(0)/\psi'(0) = \|u_0\|^2/(u_0, v_0)$. We can easily deduce from (3.9)–(3.11) that the function $\psi(t) > 0$ is convex, strictly decreasing and it approaches zero at a time $T_e < T_0$. The theorem is thus proved. $\square$

REMARKS.

(1) In the theorem, the domain $\mathcal{D}$ may be unbounded or $\mathbb{R}^d$, and the Laplacian $\nabla^2$ in equation (2.2) may be replaced by an uniformly elliptic self-adjoint operator $\sum_{i,j=1}^{d} \partial_{x_i}(a_{ij}(x)\partial_{x_j})$ with smooth coefficients $a_{ij} = a_{ji}$, and

(2) the noise term $\sigma(u, Du, x, t)\partial_t W(x, t)$ can be replaced by $\sum_{i=1}^{n} \sigma_i(u, Du, x, t)\partial_t W_i(x, t)$, where $\sigma_i(u, Du, x, t)$ are similar to $\sigma(u, Du, x, t)$, and $W_i(x, t)$ are Wiener random fields with covariance functions $EW_i(x,t)W_j(y,s) = (t \wedge s)r_{ij}(x, y)$, for $i, j = 1, \ldots, n$. In this case the expression $r(x, x)q(x, t)$ in condition (B2) should be changed to $\sum_{i=1}^{n} r_{ii}(x, x)q_i(x, t)$, where $q_i(x, t) = \sup_{\mu, \xi} \sigma_i^2(\mu, \xi, x, t)$.

**4. An example.** As an example, let us consider the following problem in the half plane $\mathcal{D} = \{x = (x_1, x_2) \in \mathbb{R}^2 | x_1 > 0\}$:

(4.1)
$$\begin{cases} \partial_t^2 u(x,t) = (c^2 \nabla^2 - \alpha)u + \lambda u^{2p-1} \\ \qquad + \sigma_0 \tan^{-1}(1 + |Du|^2)e^{-\nu t}\partial_t W(x,t), \quad t > 0, x \in \mathcal{D}, \\ u(x,0) = \dfrac{\beta}{1+|x|^2}, \quad \partial_t u(x,0) = \dfrac{1}{1+|x|^2}, \\ u(x,t)|_{x_1=0} = 0, \end{cases}$$

where $W(x, t)$ is a Wiener random field with the covariance function

(4.2) $$r(x, y) = r_0 \exp\{-\rho(x \cdot y)\} \qquad \text{for } x, y \in \mathcal{D}.$$

In the above equations, $\lambda, \beta, \rho, \nu, \sigma_0, r_0$ are all positive constants, $p \geq 2$ is an integer, and the dot product $x \cdot y = x_1 y_1 + x_2 y_2$.

We can apply Theorem 3.1 to show that the equation (4.1) has an explosive solution by choosing the above constants properly. To this end, it is easy to check that the conditions (A1)–(A3) are satisfied. It remains to verify conditions (B1)–(B3). Referring to equation (4.1), we have

(4.3) $$(u_0, v_0) = \beta \int_{\mathcal{D}} \frac{dx}{(1+|x|^2)^2} = \frac{1}{2}\beta\pi,$$

(4.4) $$(u, f(u)) = \lambda \int_{\mathcal{D}} u^{2p}\, dx,$$

(4.5) $$F(u) = \lambda \int_0^u s^{2p-1}\, ds = \frac{\lambda}{2p} u^{2p}$$

and

(4.6) $$(F(u), 1) = \frac{\lambda}{2p} \int_{\mathcal{D}} u^{2p}\, dx.$$

In view of (4.3), (4.4) and (4.6), the conditions (B1) and (B3) are clearly met. Finally, to enforce condition (B2), we must compute several integrals.

From (4.6), we get

$$(4.7) \qquad (F(u_0), 1) = \frac{\lambda}{2p} \int_{\mathcal{D}} \left(\frac{\beta}{1+|x|^2}\right)^{2p} dx = \frac{\lambda \pi \beta^{2p}}{4p(2p-1)}.$$

Since $\sigma = \sigma_0 \tan^{-1}(1 + |Du|^2)e^{-\nu t}$, we have

$$(4.8) \qquad q(x,t) = \sigma_0^2 \sup_{\xi} [\tan^{-1}(1+|\xi|^2)]^2 e^{-2\nu t} = \left(\frac{\sigma_0 \pi}{2}\right)^2 e^{-2\nu t},$$

so that

$$(4.9) \qquad \int_0^\infty \int_{\mathcal{D}} r(x,x) q(x,t)\, dx\, dt = r_0 \left(\frac{\sigma_0}{2}\pi\right)^2 \int_0^\infty \int_{\mathcal{D}} e^{-\rho|x|^2 - 2\nu t} \, dx\, dt$$
$$= \frac{r_0 \sigma_0^2 \pi^3}{16 \rho \nu}.$$

To obtain the initial energy $\mathbf{e}(u_0; v_0)$, we need to evaluate the following integrals

$$(4.10) \qquad \|v_0\|^2 = \int_{\mathcal{D}} \frac{dx}{(1+|x|^2)} = \frac{\pi}{2}, \qquad \|u_0\|^2 = \beta^2 \|v_0\|^2 = \frac{\beta^2 \pi}{2}$$

and

$$(4.11) \qquad \|Du_0\|^2 = 4\beta^2 \int_{\mathcal{D}} \frac{|x|^2}{(1+|x|^2)^4}\, dx = \frac{\pi}{3} \beta^2.$$

By the definition of the energy function (2.7), we can obtain from (4.10) and (4.11) that

$$(4.12) \qquad \mathbf{e}(u_0; v_0) = \frac{\pi}{2} \{1 + (\alpha + 2c^2/3)\beta^2\}.$$

After substituting (4.7), (4.9) and (4.12) into the condition (B2) and simplifying the result, it gives the following inequality

$$(4.13) \qquad \lambda > \frac{p(2p-1)}{\beta^{2p}} \left\{ 1 + (\alpha + 2c^2/3)\beta^2 + \frac{r_0 \sigma_0^2 \pi^2}{8\rho\nu} \right\}.$$

The above inequality shows that, for a sufficiently large amplitude $\lambda$ of the nonlinear term, keeping all other constants fixed, the mean-square solution $E\|u_t\|^2$ of equation (4.1) will explode at a finite time. Alternatively, regarding $\beta$ as a parameter and fixing all other constants, then the solution will also explode if initial displacement $u_0$ is large enough.

**Acknowledgments.** The author wishes to thank Professor J. B. Keller of Stanford University for suggesting this problem for investigation. He thanks Professors Mike Cranston and John Lowengrub for their warm hospitality during the author's visit in Fall/2008 to University of California–Irvine, where this work was completed. The author is also grateful to the referee for several helpful suggestions for improving the presentation of this paper.

# REFERENCES

[1] CHOW, P.-L. (2002). Stochastic wave equations with polynomial nonlinearity. *Ann. Appl. Probab.* **12** 361–381. MR1890069
[2] CHOW, P.-L. (2006). Asymptotics of solutions to semilinear stochastic wave equations. *Ann. Appl. Probab.* **16** 757–789. MR2244432
[3] CHOW, P.-L. (2007). *Stochastic Partial Differential Equations. Chapman & Hall/CRC Applied Mathematics and Nonlinear Science Series*. Chapman & Hall/CRC, Boca Raton, FL. MR2295103
[4] GLASSEY, R. T. (1973). Blow-up theorems for nonlinear wave equations. *Math. Z.* **132** 183–203. MR0340799
[5] JOHN, F. (1990). *Nonlinear Wave Equations, Formation of Singularities. University Lecture Series* **2**. Amer. Math. Soc., Providence, RI. MR1066694
[6] KELLER, J. B. (1957). On solutions of nonlinear wave equations. *Comm. Pure Appl. Math.* **10** 523–530. MR0096889
[7] MUELLER, C. (1997). Long time existence for the wave equation with a noise term. *Ann. Probab.* **25** 133–151. MR1428503

DEPARTMENT OF MATHEMATICS
WAYNE STATE UNIVERSITY
DETROIT, MICHIGAN 48202
USA
E-MAIL: plchow@math.wayne.edu